\documentclass[12pt]{article}

\usepackage{amssymb,amsbsy,amsthm,amsmath,graphicx,mathrsfs,epsfig,graphics,times}

\begin{document}

\title{On the geometry of a class of invariant measures and a problem of Aldous}
\author{Tim Austin}
\date{}

\maketitle


\newenvironment{nmath}{\begin{center}\begin{math}}{\end{math}\end{center}}

\newtheorem{thm}{Theorem}[section]
\newtheorem{lem}[thm]{Lemma}
\newtheorem{prop}[thm]{Proposition}
\newtheorem{cor}[thm]{Corollary}
\newtheorem{conj}[thm]{Conjecture}
\newtheorem{dfn}[thm]{Definition}
\newtheorem{prob}[thm]{Problem}


\newcommand{\A}{\mathcal{A}}
\newcommand{\B}{\mathcal{B}}
\newcommand{\C}{\mathcal{C}}
\newcommand{\E}{\mathcal{E}}
\newcommand{\calG}{\mathcal{G}}
\renewcommand{\Pr}{\mathrm{Pr}}
\newcommand{\s}{\sigma}
\newcommand{\frH}{\mathfrak{H}}
\renewcommand{\P}{\mathcal{P}}
\renewcommand{\O}{\Omega}
\newcommand{\w}{\omega}

\renewcommand{\S}{\Sigma}
\newcommand{\T}{\mathrm{T}}
\newcommand{\co}{\mathrm{co}}
\newcommand{\e}{\mathrm{e}}
\renewcommand{\d}{\mathrm{d}}
\renewcommand{\l}{\lambda}
\newcommand{\U}{\mathcal{U}}
\newcommand{\G}{\Gamma}
\newcommand{\g}{\gamma}
\newcommand{\eps}{\varepsilon}
\renewcommand{\L}{\Lambda}
\newcommand{\Sym}{\mathrm{Sym}}
\newcommand{\F}{\mathcal{F}}
\renewcommand{\a}{\alpha}
\renewcommand{\b}{\beta}
\renewcommand{\k}{\kappa}
\newcommand{\bbN}{\mathbb{N}}
\newcommand{\bbR}{\mathbb{R}}
\newcommand{\bbC}{\mathbb{C}}
\newcommand{\bbZ}{\mathbb{Z}}
\newcommand{\bbT}{\mathbb{T}}
\newcommand{\bbF}{\mathbb{F}}
\newcommand{\sfE}{\mathsf{E}}
\newcommand{\sfP}{\mathsf{P}}

\newcommand{\rsqto}{\rightsquigarrow}
\newcommand{\into}{\hookrightarrow}

\newcommand{\Atom}{\mathrm{Atom}}

\newcommand{\bb}[1]{\mathbb{#1}}
\renewcommand{\bf}[1]{\mathbf{#1}}
\renewcommand{\rm}[1]{\mathrm{#1}}
\renewcommand{\cal}[1]{\mathcal{#1}}
\newcommand{\bs}[1]{\boldsymbol{#1}}
\newcommand{\para}[1]{\paragraph{#1}}

\newcommand{\fin}{\nolinebreak\hspace{\stretch{1}}$\lhd$}
\newcommand{\tick}{\nolinebreak\hspace{\stretch{1}}$\surd$}

\parskip 7pt

\begin{abstract}
In his survey~\cite{Ald85} of notions of exchangeability, Aldous
introduced a form of exchangeability corresponding to the symmetries
of the infinite discrete cube, and asked whether these exchangeable
probability measures enjoy a representation theorem similar to those
for exchangeable sequences~\cite{HewSav55},
arrays~\cite{Hoo79,Hoo82,Ald81,Ald82} and set-indexed
families~\cite{Kal92}. In this note we to prove that, whereas the
known representation theorems for different classes of partially
exchangeable probability measure imply that the compact convex set
of such measures is a Bauer simplex (that is, its subset of extreme
points is closed), in the case of cube-exchangeability it is a copy
of the Poulsen simplex (in which the extreme points are dense). This
follows from the arguments used by Glasner and Weiss' for their
characterization in~\cite{GlaWei97} of property (T) in terms of the
geometry of the simplex of invariant measures for associated
generalized Bernoulli actions.

The emergence of this Poulsen simplex suggests that, if a
representation theorem for these processes is available at all, it
must take a very different form from the case of set-indexed
exchangeable families.
\end{abstract}

\parskip 0pt

\tableofcontents

\parskip 7pt
\parindent 0pt

\section{Introduction}

Suppose that $K$ is a standard Borel space with $\s$-algebra $\S_K$,
that $T$ is a countably infinite set and $\G$ a group of
permutations of $T$ and that $\mu$ is a probability measure on the
(standard Borel) product measurable space $(K^T,\S_K^{\otimes T})$.
Let us also always assume that $\G$ has only infinite orbits in $T$.
Then following Aldous~\cite{Ald85} we shall write that $\mu$ is
\textbf{$(T,\G)$-exchangeable} if it is invariant under the
(contravariant) coordinate-permuting action $\tau$ of $\G$ on $K^T$
given by
\[\tau^\g\big((\w_t)_{t\in T}\big) := (\w_{\g(t)})_{t\in T},\]
which is clearly measurable and invertible.  We write $\Pr^\G K^T$
for the set of all such exchangeable probability measures. We shall
sometimes refer to the index-set action $\G \curvearrowright T$ as
an \textbf{exchangeability context}.

The prototypical examples of exchangeability are arguably those of
\textbf{hypergraph exchangeability}, for which $T = \binom{S}{k}$,
the set of all $k$-subsets of a countably infinite `vertex set' $S$,
and $\G = \Sym_0(S)$, the group of all finitely-supported
permutations of $S$ acting on $T$ by vertex-permutations.  In this
case we can interpret $\mu$ as the law of a random `colouring' of
the complete $k$-uniform hypergraph on $S$ by points from the space
$K$ of `colours'.

In the simplest case $k = 1$ (so $T = S$), the precise structure of
all possible hypergraph-exchangeable measures follows from classical
theorems of de Finetti and Hewitt \& Savage (see, for
example,~\cite{HewSav55}). More recently, the case of more general
$k$ was studied by Hoover~\cite{Hoo79,Hoo82},
Aldous~\cite{Ald81,Ald82,Ald85} and Kallenberg~\cite{Kal92}, along
with a number of further extensions that are still closely related
to this hypergraph-colouring setting, leading to a more elaborate
conception of `exchangeability theory'. It turns out that in these
contexts too the exchangeable probability measures admit a
more-or-less complete structural description, albeit involving
increasingly complicated ingredients as $k$ increases: they can all
be represented as images of certain other exchangeable processes
whose laws take a particular simple form. We refer the reader
to~\cite{Aus--ERH} for a recent survey of these results and their
relations to various questions in graph and hypergraph theory, and
to the survey~\cite{Ald85} of Aldous for a general introduction to a
broader range of exchangeability contexts and to the recent book of
Kallenberg~\cite{Kal05} for the modern state of the theory.

We will not recount the details of these representation theorems
here. Rather, our interest lies in a different exchangeability
context, proposed by Aldous as a possible object of further study in
Section 16 of~\cite{Ald85}: that of \textbf{cube-exchangeability}.
Let $\bbF_2 = \{0,1\}$ be the field of two elements, and in the
$d$-dimensional vector space $\bbF_2^d$ over $\bbF_2$ write
$e_1,e_2,\ldots,e_d$ for the standard basis. Now take $T$ to be the
set $\bbF_2^{\oplus\bbN}$ of all strings of $0$s and $1$s with only
finitely many of the latter, and let $\G$ be the group of
permutations of $T$ generated by finitely-supported permutations of
the underlying copy of $\bbN$ together with all `bit-flips':
$\s_i:\bbF_2^{\oplus \bbN} \to \bbF_2^{\oplus \bbN}: x \mapsto x +
e_i$.

In this context, given any standard Borel space $K$ we shall call a
probability measure $\mu$ on $K^T$ \textbf{cube-exchangeable} if it
is invariant under the coordinate-permuting action of the above
group $\G$.  Note that we may describe this group as follows: $T$
may be written as the increasing union $\bigcup_{n\geq 1}T_n$ of the
discrete cubes $T_n := \bbF_2^n$, and now (bearing in mind our
restriction to \emph{finitely-supported} permutations of $\bbN$)
every member $g \in \G$ actually maps $T_n$ onto itself for all
sufficiently large $n$. It is easy to see that in this case a
permutation of $T_n$ is induced by a member of $\G$ if and only if
it is an isometry of $T_n$ when this latter is identified with the
$n$-dimensional Hamming cube $\{0,1\}^n$. For this reason we shall
refer to $\G$ as the group of \textbf{isometries of the
infinite-dimensional discrete cube} and denote it by
$\rm{Isom}\,\bbF_2^{\oplus \bbN}$. Note that, as in the setting of
hypergraph-exchangeability, the acting group $\G$ is locally finite
(that is, any finite collection of its elements generates a finite
subgroup); but \emph{un}like in that setting most elements of the
group (to be precise, all that involve a nontrivial translation) do
move infinitely many points of $T$.

In view of the success of the basic theory of
hypergraph-exchangeability, Aldous asked in~\cite{Ald85} whether a
similarly precise structural description is available for the class
of cube-exchangeable probability measures.  In this note we will
provide some evidence to suggest that such a structural description
may \emph{not} be available in this context --- at least not in the
very explicit form familiar from the hypergraph setting --- in the
following `soft' sense. First, we note that, provided $\G$ is
amenable (as it certainly is in our examples), the basic
representation theorems for hypergraph exchangeable laws fall into a
certain quite general pattern, and that this pattern has, in
particular, the consequence that for a compact metric $K$ the set of
all extreme points (that is, ergodic members) of $\Pr^\G K^T$ forms
a closed subgroup of this compact convex set in the vague topology;
that is, this convex set is a \textbf{Bauer simplex}. On the other
hand, we will show that provided $K$ is not a singleton, this set
$\Pr^\G K^T$ in the case of cube-exchangeability has the very
different property of being a copy of the \textbf{Poulsen simplex}:
its extreme points form a vaguely dense subset.  This suggests that
any representation theorem describing this set, if one is available,
must take a rather different form from the earlier set-indexed
examples.

\subsubsection*{Remark on notation}

Our basic combinatorial and measure-theoretic notation is completely
standard.  If $(X,\rho)$ is a metric space, $x,y \in X$ and $\eps >
0$, we shall sometimes write $x\approx_\eps y$ in place of
$\rho(x,y) < \eps$ when the particular metric $\rho$ is understood.

\subsubsection*{Acknowledgements}

My thanks go to Terence Tao and Yehuda Shalom for helpful
discussions.

\section{The form of previous representation theorems for exchangeable
measures}\label{sec:repthms}

In this section we introduce a general template for a kind of
representation theorem for exchangeable laws, which in particular
characterizes the basic representation theorems for the cluster of
variations on hypergraph-exchangeability.

These theorems all focus on representing an arbitrary
$(T,\G)$-exchangeable process as an image (in a suitable sense) of
another exchangeable process (possibly with a different index set)
for which the different random variables are all mutually
independent.

\begin{dfn}[Ingredients]\label{dfn:ingredients}
Let $\G \curvearrowright T$ be an exchangeability context and $K$ a
fixed compact metric space.  By a list of \textbf{representation
data} we understand:
\begin{itemize}
\item a sequence of auxiliary index sets $T_1$, $T_2$, \ldots each endowed with some action
$\G \curvearrowright T_i$ that has \emph{only infinite orbits};
\item a disjoint sequence of \textbf{dependency maps} $\phi_i:T\to
\binom{T_i}{<\infty}$ that are $\G$-covariant, in that
$\phi_i(\g(t)) = \g(\phi_i(t))$;
\item and a family of probability kernels
\[\k_t:[0,1] \times [0,1]^{\phi_1(t)}\times [0,1]^{\phi_2(t)} \times
\cdots \rightsquigarrow K\] that is $\G$-covariant, in that
$\k_{\g(t)} = \k_t\circ(\rm{id}_{[0,1]}\times \tau_1^\g\times
\cdots)$.
\end{itemize}

Given ingredients as above, we denote by $\k^{(T)}$ the kernel
$[0,1]\times [0,1]^{T_1}\times [0,1]^{T_2}\times \cdots
\rightsquigarrow K^T$ given by
\[\k^{(T)}(x_0,\bf{x}_1,\ldots,\,\cdot\,) = \bigotimes_{t\in T}\k_t(x_0,\bf{x}_1|_{\phi_1(t)},\ldots,\,\cdot\,).\]
\end{dfn}

Of the conditions on the data introduced above, perhaps the least
intuitive is that the actions $\G \curvearrowright T_i$ may not have
finite orbits (although it certainly holds in the case of hypergraph
exchangeability); we shall later need to play this off against the
finiteness of the sets $\phi_i(t)$, and it does hold for the case of
hypergraph-exchangeability.

Now and henceforth we will denote by $\mu_{\rm{L}}$ Lebesgue measure
on the unit interval $[0,1]$, and by the shorthand
$\mu_\rm{L}^{\ast\otimes T_1\otimes T_2\otimes \cdots}$ the product
measure $\mu_\rm{L}\otimes\mu_\rm{L}^{\otimes T_1}\otimes
\mu_\rm{L}^{\otimes T_2}\otimes\cdots$.

\begin{dfn}[Representability]\label{dfn:representability}
Given an exchangeability context $\G \curvearrowright T$ and a
compact metric space $K$, we shall say that a $(T,\G)$-exchangeable
law $\mu \in \Pr\,K^T$ is \textbf{representable} if there is a list
of ingredients as above, with only the kernels $\k_t$ allowed to
depend on $\mu$ or $K$, such that $\mu =
\k^{(T)}_\#(\mu_\rm{L}\otimes\mu_\rm{L}^{\otimes T_1}\otimes
\mu_\rm{L}^{\otimes T_2}\otimes\cdots)$.

If an exchangeability context $(T,\G)$ is such that all exchangeable
laws on $K^T$ are representable for any compact metric $K$ then we
shall say that $(T,\G)$ \textbf{always admits representation}.
\end{dfn}

We must stress that our chosen definition of representability is not
completely canonical: although we are guided by the classical
representation theorems for hypergraph-exchangeable laws and their
relatives, these leading examples are sufficiently closely related
one to another that it is not quite clear which features of their
representation theorems we should try to keep, and which to discard,
when abstracting to a more general definition.  The choice we have
made seems to be simple and natural, and also to reflect many of the
uses to which these representation theorems are put
(see~\cite{Kal05}), but certainly it has also been selected partly
because it works for what follows. An alternative formulation of the
representation theorem for exchangeable arrays can be given instead
in terms, for example, of sequences of auxiliary compact metric
spaces $Z_0,Z_1,Z_2,\ldots$ and index sets $T_1,T_2,\ldots$ with
$\G$-actions $\a_1,\a_2,\ldots$ from which all exchangeable laws are
then obtained as pushforwards of probability measures on the product
space $Z_0 \times Z_1^{T_1}\times \cdots$ that are invariant under
the associated overall coordinate-permuting action of $\G$ and have
the additional property that the coordinates in $Z_{i+1}$ are
conditionally independent given the coordinates in every $Z_j$ for
$j \leq i$. The representation theorem for exchangeable arrays is
treated in these terms, for example, in~\cite{Aus--ERH}, where this
choice is dictated by the use to which that theorem is then put in
Section 3 of~\cite{AusTao--hereditarytest}; however, the formalism
of representability extracted this way seems much less amenable to
our needs, as well as further from the classical descriptions of
Aldous and Kallenberg, and so we have settled for the above instead.

In our present terms the main Representation Theorem of Aldous,
Hoover and Kallenberg for hypergraph-exchangeable laws with $T :=
\binom{S}{k}$ and $\G := \Sym_0(S)$ with its canonical action may be
written as follows.

\begin{thm}[Representation Theorem for hypergraph-exchangeable laws]
A hypergraph-exchangeable law $\mu$ is representable using the data
$T_i := \binom{S}{i}$ for $i \leq k$ and $T_{k+1} = T_{k+2} = \ldots
= \rm{triv}.$, the dependency maps $\phi_i:t \mapsto \binom{t}{i}$
for $t \in \binom{S}{k}$ and $i \leq k$ and $\phi_i \equiv
\emptyset$ if $i \geq k+1$, and some deterministic maps $\k_t$ that
depend on the particular choice of $\mu$. \qed
\end{thm}

Although we have allowed arbitrary probability kernels $\k_t$ in our
present formalism, in the above concrete representation theorem (and
its relatives in such works as~\cite{Ald82,Kal92}) they are all
deterministic maps. However, a simple transfer argument shows that
this difference is purely cosmetic.

\begin{lem}
A $(T,\G)$-exchangeable law is representable if and only if it is
representable using deterministic maps $\k_t:[0,1] \times
[0,1]^{\phi_1(t)}\times [0,1]^{\phi_2(t)} \times \cdots \to K$.
\end{lem}

\textbf{Proof}\quad Clearly representability using deterministic
maps amounts to a special case of representability, so we need only
prove that any representable law is representable using
deterministic maps.  However, if we have a list of ingredients that
represents $\mu$ with kernels $\k_t$, then by the standard Transfer
Theorem (Theorem 6.10 in Kallenberg~\cite{Kal02}) we may find
deterministic maps
\[\theta_t:[0,1]\times \big([0,1] \times
[0,1]^{\phi_1(t)}\times [0,1]^{\phi_2(t)} \times \cdots\big) \to K\]
such that \[\k_t(x_0,\bf{x}_1,\ldots,\,\cdot\,) = \mu_\rm{L}\big\{y
\in [0,1]:\ \theta_t(y,x_0,\bf{x}_1,\ldots) \in\,\cdot\,\big\}.\]
Now, as is standard, the Lebesgue spaces $([0,1],\mu_\rm{L})$ and
$([0,1]\times [0,1],\mu_\rm{L}\otimes\mu_\rm{L})$ are isomorphic,
say via the Borel map $\xi:[0,1]\to [0,1]^2$, and so now defining
\[\tilde{\k}_t(x_0,\bf{x}_1,\ldots) := \theta_t(\xi(x_0),\bf{x}_1,\ldots)\]
we can check at once from the above relations that these
deterministic maps also represent the original law $\mu$. \qed

\section{Bauer simplices from exchangeability}

We will now prove that if $\G$ is amenable and the exchangeability
context $(T,\G)$ always admits representation then its simplices
$\Pr^\G K^T$ of exchangeable laws must be Bauer for any $K$. We will
also give a direct deduction of this Bauer property in the
representative example of hypergraph exchangeability without using
representability, both for completeness and because it seems
interesting to compare this direct proof with arguments to prove the
Poulsen property in the case of cube-exchangeability in the next
section.

\subsection{The Bauer property from representability}

\begin{lem}
If $\G$ is amenable, and if an $(T,\G)$-exchangeable probability
measure $\mu \in \Pr^\G K^T$ is representable at all, then it is
ergodic if and only if it is representable by kernels $\k_t$ not
depending on the first coordinate.
\end{lem}

\textbf{Proof}\quad First suppose that $\mu$ is ergodic, and write
it as $\k^{(T)}_\#\mu_\rm{L}^{\ast\otimes T_1\otimes
T_2\otimes\cdots}$ for some suitable family $\k$.  Now define the
family $\k^u$ of kernels $\k^u_t:[0,1]\times [0,1]^{\phi_1(t)}\times
[0,1]^{\phi_2(t)}\times \cdots \rightsquigarrow K$ by
$\k^u_t(x_0,\bf{x}_1,\ldots) := \k_t(u,\bf{x}_1,\ldots)$ (this makes
sense and is unambiguous up to equality for almost every $u$);
clearly none of these depends on the first coordinate in
$[0,1]\times [0,1]^{\phi_1(t)}\times\cdots$, and also each
$(\k^u)^{(T)}_\#\mu_\rm{L}^{\ast\otimes T_1\otimes T_2\otimes
\cdots}$ is another $\G$-invariant probability on $K^T$ such that
\[\mu = \int_0^1(\k^u)^{(T)}_\#\mu_\rm{L}^{\ast\otimes T_1\otimes T_2\otimes
\cdots}\,\d u.\] By the ergodicity of $\mu$ this decomposition must
be trivial, and so $(\k^u)^{(T)}_\#\mu_\rm{L}^{\ast\otimes
T_1\otimes T_2\otimes \cdots} = \mu$ for almost-every $u$; hence
almost any of the kernel families $\k^u$ will suffice.

Now suppose, on the other hand, that each $\k_t$ does not depend on
the first coordinate in $[0,1]\times [0,1]^{T_1}\times \cdots$, and
that $A,B \subseteq K^T$ are two Borel finite-dimensional cylinder
sets, say determined by the finite sets of coordinates $I,J$
respectively.  Then by our assumption that all orbits of $\G$ on $T$
and on $T_i$ are infinite and that $\G$ is amenable, it follows that
for some density-$1$ subset of $F \subset \G$ we have
$\phi_i(g(I))\cap \phi_i(J) = \emptyset$ for all $g \in F$. However,
this implies that $\k_t$ and $\k_s$ have no arguments in common for
$t \in g(I)$ and $s \in J$, and so the sets $\tau^g(A)$ and $B$ must
be independent under $\mu$.  In fact this proves not only
ergodicity, but even weak mixing, and we are done. \qed

\begin{prop}[Representability implies
Bauer]\label{prop:Bauerfromrecipe} If $\G$ is amenable and the
exchangeability context $(T,\G)$ always admits representation then
$\Pr^\G K^T$ is a Bauer simplex for any compact metric $K$.
\end{prop}

\textbf{Proof}\quad We know that $\Pr^\G K^T$ is a compact convex
set and that its extreme points are precisely those members that can
be represented by some collection of kernels $\k_t$ not depending on
the first coordinate in $[0,1]\times [0,1]^{\phi_1(t)}\times
\cdots$; thus we need only show that if $\mu_n =
(\k_n)^{(T)}_\#\mu_\rm{L}^{\ast\otimes T_1\otimes T_2\otimes
\cdots}$ are a vaguely convergent sequence of such measures then
their limit $\mu$ admits a similar representation.

However, for each $t$ the kernel $\k_t$ defines a joining of the
probability measures $\mu_\rm{L}^{\ast\otimes T_1\otimes T_2\otimes
\cdots}$ and $(\pi_t)_\#\mu_n$ on the product space $[0,1]\times
[0,1]^{\phi_1(t)}\times \cdots \times K$ under which the very first
coordinate is independent from all the others (because $\k_t$ does
not depend on this coordinate), and so, passing to a subsequence if
necessary, we may assume that these joinings also converge to some
fixed probability measure $\l_{\infty,t}$ on this product space. It
is clear that this new measure will still have projection onto
$[0,1]\times [0,1]^{\phi_1(t)}\times\cdots$ equal to
$\mu_\rm{L}^{\ast\otimes T_1\otimes T_2\otimes \cdots}$ and will
still enjoy the independence of the first coordinate from everything
else, and so if we now disintegrate these $\l_{\infty,t}$ over that
first projection we recover kernels $\k_{\infty,t}$ that also do not
depend on the very first coordinate and represent $\mu$, as
required. \qed

\textbf{Remark}\quad I do not know whether the assumption of
amenability could be removed from the preceding arguments. \fin

\subsection{The Bauer property in the particular context of hypergraph exchangeability}

Before moving on, let us include a second proof that the classical
hypergraph-exchangeability context has the Bauer property that uses
only a very elementary property enjoyed by that context, rather than
the representation theorem. This subsection is not essential to the
main thread of this note, but is included mainly to advertise the
question of whether the argument that it contains can be generalized
further.

\begin{dfn}[Distant multiple transitivity]
We shall write that an exchangeability context $(T,\G)$ is
\textbf{distantly multiply transitive (DMT)} if for any finite $I,J
\subset T$ there is some subset $E \subseteq \G$ of density $1$ and
such that for any $\g_1,\g_2 \in E$ there is some $\xi \in \G$ with
$\xi\!\upharpoonright_I = \rm{id}_I$ and $\xi\circ
\g_1\!\upharpoonright_J = \g_2\!\upharpoonright_J$.
\end{dfn}

It is immediate to check that the hypergraph exchangeability context
is DMT, and so the following result applies to that context in
particular.

\begin{prop}[DMT implies Bauer]
If $\G$ is amenable and $(T,\G)$ is DMT then it has the Bauer
property.
\end{prop}

\textbf{Proof}\quad We follow closely the analogous argument of
Glasner and Weiss in~\cite{GlaWei97}. Suppose that $\G$ is amenable,
that $(T,\G)$ is DMT, that $\mu \in \Pr^\G K^T$ can be vaguely
approximated by ergodic measures, and that $A \in \S_{K^T}$ is
invariant with $a := \mu(A) \in [0,1]$. For any $\eps
> 0$ there are a finite set $J \subset T$ and a continuous function
$f:K^J \to [0,1]$ such that $\|1_A - f\circ\pi_J\|_{L^1(\mu)} <
\eps$, and hence $\int_{K^T}f\circ\pi_J\,\d\mu\approx_\eps a$. From
the invariance of $A$ it follows that we actually have $\|1_A -
f\circ\pi_J\circ\tau^\g\|_{L^1(\mu)} < \eps$ for any $\g \in \G$.

Now, since $(T,\G)$ is DMT and $J$ is finite, there is some $E
\subseteq \G$ with asymptotic density $1$ such that for any
$\g_1,\g_2 \in E$ there is some $\xi \in \G$ such that
$\xi\!\upharpoonright_J = \rm{id}_J$, and so
$f\circ\pi_J\circ\tau^\xi = f\circ\pi_J$, whereas
$\xi\circ\g_1\!\upharpoonright_J = \g_2\!\upharpoonright_J$ and so
$f\circ\pi_J\circ\tau^{\g_1}\circ\tau^\xi =
f\circ\pi_J\circ\tau^{\g_2}$. Let us now fix some representative
member $\g_0 \in E$.

Next, since $f\circ\pi_J$ and
$(f\circ\pi_J)\cdot(f\circ\pi_J\circ\tau^{\g_0})$ are continuous, by
assumption we can always find some ergodic $\mu' \in \Pr^\G K^T$
with
\[\int_{K^T}f\circ\pi_J\,\d\mu' \approx_\eps
\int_{K^T}f\circ\pi_J\,\d\mu\] and
\[\int_{K^T}(f\circ\pi_J)\cdot(f\circ\pi_J\circ\tau^{\g_0})\,\d\mu'
\approx_\eps
\int_{K^T}(f\circ\pi_J)\cdot(f\circ\pi_J\circ\tau^{\g_0})\,\d\mu.\]

Letting $(I_n)_{n \geq 1}$ be a F\o lner sequence in $\G$, it
follows from the ergodicity of $\mu'$ that
\[\frac{1}{|I_N|}\sum_{\g \in I_N}\int_{K^T}(f\circ\pi_J)\cdot(f\circ\pi_J\circ\tau^\g)\,\d\mu' \to \Big(\int_{K^T}f\circ\pi_J\,\d\mu'\Big)^2 \approx_{2\eps} \Big(\int_{K^T}f\circ\pi_J\,\d\mu\Big)^2\approx_{2\eps} a^2\]
as $N \to \infty$. On the other hand, we know that for $N$
sufficiently large at least $(1-\eps)$-proportion of $\g \in I_N$
lie in $E \cap I_N$, and that $\g_0 \in E\cap I_N$, and so by
choosing a suitable $\xi$ they must all give exactly the same value
for $\int_{K^T}(f\circ\pi_J)\cdot(f\circ\pi_J\circ\tau^\g)\,\d\mu'$;
and therefore for $N$ sufficiently large we must also have
\begin{multline*}
\frac{1}{|I_N|}\sum_{\g \in
I_N}\int_{K^T}(f\circ\pi_J)\cdot(f\circ\pi_J\circ\tau^\g)\,\d\mu'
\approx_\eps
\int_{K^T}(f\circ\pi_J)\cdot(f\circ\pi_J\circ\tau^{\g_0})\,\d\mu'\\
\approx_\eps
\int_{K^T}(f\circ\pi_J)\cdot(f\circ\pi_J\circ\tau^{\g_0})\,\d\mu\approx_{2\eps}
\int_{K^T}1_A\cdot 1_A\,\d\mu = \mu A = a.
\end{multline*}

Combining these approximations shows that $a \approx_{6\eps} a^2$
for any $\eps > 0$, and so in fact we must have $a \in \{0,1\}$, and
$\mu$ must itself be ergodic. \qed

\section{The Poulsen property for cube-exchangeable measures}

We will now show that, quite unlike the cases studied in the
previous two sections, if $(T,\G)$ is the cube-exchangeability
context (and $K$ is nontrivial) then $\Pr^\G K^T$ is actually the
Poulsen simplex. This argument is also closely motivated by that of
Glasner and Weiss in~\cite{GlaWei97}, where they show that in the
case of the exchangeability context $(\G,R_\G)$ comprising a group
$\G$ and its right-regular representation on itself, the simplex
$\Pr^\G \{0,1\}^\G$ of invariant probability measures is either
Bauer or Poulsen precisely according as $\G$ has or fails Kazhdan's
property (T).  No condition like property (T) will enter our
analysis --- indeed, the groups of immediate interest to us are all
locally finite, hence trivially amenable --- but we will follow
closely the basic steps of their construction.

There are essentially two of these steps.  We first show that in
case $K = \{0,1\}$ the particular example $\frac{1}{2}\delta_\bf{0}
+ \frac{1}{2}\delta_\bf{1}$ of a non-ergodic member of $\Pr^\G
\{0,1\}^T$ is vaguely approximable by members that are not only
ergodic, but actually weakly mixing; and then we use this fact
through the construction of a certain joining to show that quite
generally whenever $\mu_1$ and $\mu_2$ in $\Pr^\G K^T$ are
approximable by ergodic measures, so is their average
$\frac{1}{2}\mu_1 + \frac{1}{2}\mu_2$.  We need to ensure weak
mixing in the first step because we shall need to ensure the
ergodicity of a certain product in the second, but this makes little
difference to the other details of the proofs. It is easy to see
that this then implies the Poulsen property.

\begin{lem}\label{lem:basicvagueapprox}
Let $(T,\G) = (\bbF_2^{\oplus\bbN},\rm{Isom}\,\bbF_2^{\oplus\bbN})$.
Then the measure $\frac{1}{2}\delta_\bf{0} +
\frac{1}{2}\delta_{\bf{1}} \in \Pr^\G\{0,1\}^T$ is vaguely
approximable by weakly mixing members of $\Pr^\G\{0,1\}^T$.
\end{lem}

\textbf{Proof}\quad We need to show that for any $\eps > 0$ and $N
\geq 1$ there is some strongly mixing measure $\mu \in
\Pr^\G\{0,1\}^T$ such that both
\[\mu\{\w \in \{0,1\}^T:\ \w\!\upharpoonright_{\bbF_2^N} = \bf{0}\} \geq \frac{1}{2} - \eps\]
and
\[\mu\{\w \in \{0,1\}^T:\
\w\!\upharpoonright_{\bbF_2^N} = \bf{1}\} \geq \frac{1}{2} - \eps.\]

There are many possible ways to construct such a $\mu$; the
following seems to be one of the simplest.  We specify $\mu$ as the
law of the member of $\{0,1\}^T$ output by the following random
procedure.  For any $p \in [0,1]$ let $\nu_p$ be the product measure
on $\bbF_2^\bbN$ with $\nu_p\{z:\ z_i = 1\} = p$ for every $i \in
\bbN$; and for any $z = (z_i)_{i \in \bbN}\in \bbF_2^\bbN$ and $x
\in \bbF_2^{\oplus \bbN}$ define $\langle x,z\rangle := \sum_{i \in
\bbN}x_iz_i \mod 2$ (this sum being actually always finite). Now let
$\mu$ be the law of the characteristic function of the random subset
$\{x\in \bbF_2^{\oplus \bbN}:\ \langle x,z\rangle + \eta = 0 \mod
2\}$ where $z \sim \nu_p$ for some very small $p
> 0$ and $\eta \in \bbF_2$ is chosen independently and uniformly at
random.

It is clear that this $\mu$ is $\G$-invariant and strongly mixing
provided $p \neq 0$, but if $p$ is very small then for our chosen
$N$ we have $\nu_p\{z:\ z_1 = z_2 = \ldots = z_N = 0\} \geq 1 -
\eps$, and conditioned on the event $\{z:\ z_1 = z_2 = \ldots = z_N
= 0\}$ we must have also
\[1_{\{x\in \bbF_2^N:\ \langle x,z\rangle +
\eta = 0 \mod 2\}} = \left\{\begin{array}{ll} \bf{1}&\quad\hbox{if
}\eta =
0\quad\hbox{(occurs with prob. $\frac{1}{2}$)}\\
\bf{0}&\quad\hbox{if }\eta = 1\quad\hbox{(occurs with prob.
$\frac{1}{2}$)},\end{array}\right.\] which proves the desired vague
approximation to $\frac{1}{2}\delta_\bf{0} +
\frac{1}{2}\delta_\bf{1}$. \qed

\begin{thm}\label{thm:Poulsen}
The cube-exchangeability context $(T,\G) = (\bbF_2^{\oplus
\bbN},\rm{Isom}\,\bbF_2^{\oplus\bbN})$ has the Poulsen property.
\end{thm}

\textbf{Proof}\quad Let $K$ be any compact metric space containing
at least two points. As argued by Glasner and Weiss
in~\cite{GlaWei97}, it suffices to prove that for any two ergodic
$\mu_1, \mu_2 \in \Pr^\G K^T$, their average $\frac{1}{2}\mu_1 +
\frac{1}{2}\mu_2$ can be approximated by ergodic members of $\Pr^\G
K^T$; for then it follows by repeated approximation that the ergodic
probability measures must be dense in their own convex hull, but
this is the whole of $\Pr^\G K^T$.

Thus, it is enough to show that for any $\eps > 0$ and finite list
of continuous functions $f_1$, $f_2$, \ldots, $f_m:K^T \to [0,1]$
there is some ergodic $\mu \in \Pr^\G K^T$ such that
\[\int_{K^T}f_i\,\d\mu \approx_{2\eps} \frac{1}{2}\int_{K^T}f_i\,\d\mu_1 + \frac{1}{2}\int_{K^T}f_i\,\d\mu_2\quad\quad\forall i\leq m.\]
Moreover, by the Stone-Weierstrass Theorem we may assume each $f_i$
depends only on coordinates in some fixed finite subset $J \subset
T$, and so may factorize and rewrite it as $f_i\circ\pi_J$.

First, let us choose $\mu_0 \in \Pr^\G \{0,1\}^T$ weakly mixing and
satisfying $\mu_0 (A) \approx_\eps \frac{1}{2}\delta_\bf{0} (A) +
\frac{1}{2}\delta_\bf{1} (A)$ for all $A \subseteq \{0,1\}^T$
depending only on coordinates in $J$; this is possible by
Lemma~\ref{lem:basicvagueapprox}.  Now consider any ergodic
cube-exchangeable joining $\l$ of the two measures $\mu_1$ and
$\mu_2$ on the product space $(K^2)^T$ (such can be obtained, for
example, by taking any ergodic component of the simple product
$\mu_1\otimes \mu_2$), and now from this construct the product
$\mu_0\otimes \l$, a member of $\Pr^\G(\{0,1\}\times K^2)^T$.  Since
$\mu_0$ is weakly mixing, this product is still ergodic.

We now complete the proof by specifying a $\G$-equivariant map
$\psi:(\{0,1\}\times K^2)^T \to K^T$ whose law as a $K^T$-valued
random variable under $\mu_0\otimes \l$ will be the ergodic
approximating measure that we seek: given a point
$(\eta,\w^{(1)},\w^{(2)}) \in (\{0,1\}\times K^2)^T$, we define
$\psi(\eta,\w^{(1)},\w^{(2)})_t$ to be $\w^{(1)}_t$ if $\eta_t = 0$,
and $\w^{(2)}_t$ if $\eta_t = 1$.  Let us also write $\psi^{(1)}$
and $\psi^{(2)}$ for the usual projection maps $(\{0,1\}\times
K^2)^T \to K^T$ onto the first and second copies of $K^T$
respectively.

It is clear that this $\psi$ is equivariant, and that its law
$\psi_\#(\mu_0\otimes \l)$ must, like $\mu_0\otimes \l$, be ergodic.
Finally,
\begin{eqnarray*}
&&\int_{K^T}f_i\circ\pi_J\,\d\psi_\#(\mu_0\otimes \l) =
\int_{(\{0,1\}\times K^2)^T}f_i\circ\pi_J\circ\psi\,\d(\mu_0\otimes \l)\\
&&= \int_{\{\eta\!\upharpoonright_J =
\bs{0}\}}f_i\circ\pi_J\circ\psi\,\d(\mu_0\otimes \l) +
\int_{\{\eta\!\upharpoonright_J =
\bs{1}\}}f_i\circ\pi_J\circ\psi\,\d(\mu_0\otimes \l)\\
&&\quad + \int_{\{\eta\!\upharpoonright_J =
\bs{0}\}^\complement\cap\{\eta\!\upharpoonright_J =
\bs{1}\}^\complement}f_i\circ\pi_J\circ\psi\,\d(\mu_0\otimes \l)\\
&&\approx_\eps \mu_0\{\eta\!\upharpoonright_J =
\bs{0}\}\cdot\int_{(K^2)^T}f_i\circ\pi_J\circ\psi^{(1)}\,\d\l +
\mu_0\{\eta\!\upharpoonright_J =
\bs{1}\}\cdot\int_{(K^2)^T}f_i\circ\pi_J\circ\psi^{(2)}\,\d\l\\
&&\approx_\eps \frac{1}{2}\int_{K^T}f_i\circ\pi_J\,\d\mu_1 +
\frac{1}{2}\int_{K^T}f_i\circ\pi_J\,\d\mu_2,
\end{eqnarray*}
where we have deduced from the known quality of our approximation
$\mu_0 \approx \frac{1}{2}\delta_{\bs{0}} +
\frac{1}{2}\delta_{\bs{1}}$ that
\[\mu_0\{\eta\!\upharpoonright_J =
\bs{0}\},\ \mu_0\{\eta\!\upharpoonright_J = \bs{1}\} \approx_\eps
\frac{1}{2}\] and
\[\mu_0(\{\eta\!\upharpoonright_J =
\bs{0}\}^\complement\cap\{\eta\!\upharpoonright_J =
\bs{1}\}^\complement)\approx_\eps  0.\] This completes the proof.
\qed

\begin{cor}[Failure of cube-exchangeable representability]
For the infinite discrete cube context $(T,\G)$, the exchangeable
laws $\Pr^\G [0,1]^T$ do not admit representation.
\end{cor}

\textbf{Proof}\quad This follows at once from
Theorem~\ref{thm:Poulsen} and
Proposition~\ref{prop:Bauerfromrecipe}. \qed

\section{Some further questions}

\subsection{Further analysis of cube-exchangeable measures}

In~\cite{Ald85} (Examples 16.7 and 16.10) Aldous introduces an
interesting family of examples of cube-exchangeable measures built
from reversible random walks on a compact Abelian group, and asks
whether these might play a r\^ole in a more complete representation
theorem for such measures. Since they do not seem to fall easily
into the framework set up in Section~\ref{sec:repthms}, it would be
remiss of us not to mention them separately.

Letting $U$ be such a group endowed with its Borel $\s$-algebra
$\S_U$ and Haar measure $\mu_U$, and suppose also that $\nu \in
\Pr\,U$.  From this data we can define a measure $\mu$ on
$U^{\bbF_2^{\oplus \bbN}}$ as the law of the following randomized
selection of a point $(g_v)_{v\in \bbF_2^{\oplus \bbN}}$ of this
space:
\begin{itemize}
\item First select $g_\bf{0} \in U$ uniformly at random;
\item Now select $g^\circ_i \in U$ for each $i \in \bbN$ independently
at random with law $\nu$, and let $g_v := g_\bf{0} +
\sum_{i\in\bbN}v_ig^\circ_i$ for all $v = (v_i)_{i\in\bbN} \in
\bbF_2^{\oplus \bbN}$.
\end{itemize}

The $\Sym_0(\bbN)$-symmetry (`hypergraph-exchangeability') of this
law $\mu$ is manifest; in order to guarantee full
cube-exchangeability it turns out to be necessary and sufficient
that $\nu$ satisfy the symmetry condition that the two maps
$(g_0,g_1) \mapsto (g_0,g_0+g_1)$ and $(g_0,g_1) \mapsto
(g_0+g_1,g_0)$ have the same law under the product measure
$\mu_U\otimes \nu_0$.

Notice that we have already met one of these Abelian group examples
in the form of the measure $\mu$ constructed from $\nu_p$ during the
proof of Lemma~\ref{lem:basicvagueapprox}.

Cube-exchangeable systems of this form (or, more generally, factors
of such systems) are surely rather special, but they fit into a
considerably more general framework, and this may afford some
greater purchase over the general case. Let us approach this
generalization from a rather different direction.

Since $U$ is an Abelian group we may describe a general point of
$U^{\bbF_2^{\oplus \bbN}}$ using a M\"obius inversion formula: for
any $(g_v)_{v\in\bbF_2^{\oplus \bbN}} \in U^{\bbF_2^{\oplus \bbN}}$
there are unique $(u_\a)_{\a\in\binom{\bbN}{<\infty}} \in
U^{\binom{\bbN}{<\infty}}$ such that
\[g_v = \sum_{\a \subseteq v^{-1}\{1\}}u_\a= \sum_{\a \in \binom{\bbN}{<\infty}}\Big(\prod_{i\in\a}v_i\Big)u_\a\quad\quad\forall v\in\bbF_2^{\oplus \bbN},\]
and it is routine to check that the resulting bijection
$\Phi:U^{\bbF_2^{\oplus\bbN}}\to U^{\binom{\bbN}{<\infty}}$ is
actually a homeomorphism, and that it is covariant for the
coordinate-permuting actions of $\Sym_0(\bbN)$ on the domain and on
the target.  It follows that any hypergraph-exchangeable $\mu \in
\Pr^{\Sym_0(\bbN)} U^T$ is pushed forward by $\Phi$ to a
hypergraph-exchangeable measure $\Phi_\#\mu$ on
$U^{\binom{\bbN}{<\infty}}$, and indeed that this gives an affine
homeomorphism between the simplices of hypergraph-exchangeable
measures.  However, the stronger assumption that $\mu$ be
cube-exchangeable is then converted under $\Phi$ into a rather
larger set of additional symmetries for $\Phi_\#\mu$, and these are
not obviously easier to describe explicitly than the original
cube-exchangeable structure of $\mu$.

Indeed, if $\mu \in \Pr^\G K^T$ for an arbitrary compact metric
space $K$ and the one-dimensional marginals $(\pi_v)_\#\mu \in
\Pr\,K$ (which must all agree) are atomless, then we can simply
choose any non-discrete compact Abelian group $U$ and a function
$(K,(\pi_\bf{0})_\#\mu) \to (U,\mu_U)$ that defines a
measure-algebra-isomorphism and observe that applying this function
pointwise gives an isomorphism from $\G \curvearrowright
K^{\bbF_2^{\oplus\bbN}}$ to $\G \curvearrowright U^{\bbF_2^{\oplus
\bbN}}$, and so without any additional assumptions the above
examples of cube-exchangeable laws on Abelian groups lose no
generality at all. However, we might ask whether we can find a route
to a more interesting representation theorem through a canny choice
of the isomorphism $(K,(\pi_\bf{0})_\#\mu) \to (U,\mu_U)$, for which
the additional constraints on the joint law of
$(u_\a)_{\a\in\binom{\bbN}{< \infty}}$ can then be described
explicitly.  A little more generally, can we some $U$ and some
cube-exchangeable measure $\theta$ on $U^{\bbF_2^{\oplus\bbN}}$ of
an especially simple form such that $\mu$ is a coordinatewise factor
of $\theta$, say $\mu = (f^{\bbF_2^{\oplus\bbN}})_\#\theta$ for some
Borel $f:U \to K$.  For example, can we choose a $\theta$ under
which the summands in the M\"obius inversion formula corresponding
to sets of different sizes are independent?

We will not offer so much here, but merely note that more can be
said in certain simple cases. For example, if $u_\a = 0$ a.s.
whenever $|\a| \geq 2$, then the above laws $\mu$ must be measures
of the kind described in Aldous' example, as may be checked by hand
from the rank-$2$ case of the hypergraph-exchangeability
representation theorem applied to $\Phi_\#\mu$.

More generally, we can focus attention on the sub-simplices of
cube-exchangeable laws that are concentrated on certain
$\G$-invariant closed subsets of $K^T$. For each $r \geq 1$ let
$\O_r$ be the subset of those $g \in U^{\bbF_2^{\oplus \bbN}}$ with
the property that `all $r$-faces sum to zero':
\[g \in \O_r\quad\quad\Leftrightarrow\quad\quad \sum_{v\in F}g_v = 0\quad\hbox{for each $r$-face }F \subseteq \bbF_2^{\oplus \bbN}.\]

This suggestion is made by Aldous in~\cite{Ald85} (example 16.20),
where he also points out that some such restricted measures already
defeat any overly-simple approach to a representation theorem for
cube-exchangeability using group random walks.

It is easy to check that concentration on $\O_2$ is equivalent to
the abovementioned condition that $u_\a = 0$ a.s. whenever $|\a|
\geq 2$. It turns out that in the special case $U = \bbF_2$ the
points of $\O_r$ have a particularly simple explicit description: in
this case, identifying $U^{\bbF_2^{\oplus \bbN}}$ as the space of
functions $\bbF_2^{\oplus \bbN} \to \bbF_2$, an explicit calculation
of the M\"obius inversion gives at once that a function
$g:\bbF_2^{\oplus \bbN} \to \bbF_2$ lies in $\O_r$ if and only if it
is a polynomial of degree at most $r$. (Note that for a general
field $\bb{K}$ it is fairly straightforward to prove that those
functions $f:\bb{K}^d \to \bb{K}$ that have zero sum across any
\emph{affine} copy of the $r$-dimensional discrete cube in
$\bb{K}^d$ must be a polynomial of degree at most $r$, for any
underlying field $\bb{K}$.  However, under the present weaker
assumption of zero-sums across only \emph{isometric} copies of the
$r$-cube in $\bbF_2^d$, and it is not hard to find examples showing
that the implication of degree-$r$ polynomiality follows only over
the smallest field $\bbF_2$.)

\subsection{The geometry of subsimplices and relations to property testing}

Theorem~\ref{thm:Poulsen} has consequences for the relations between
the vague topology and the `$\bar{\rm{d}}$'- (or joining) topology
(considered by Aldous in the case of hypergraph exchangeability
in~\cite{Ald82.2} and Section 15 of~\cite{Ald85}).  This latter is
defined by the \textbf{$\bar{\rm{d}}$-metric} $\rho$ on exchangeable
probability measures, given by
\[\rho(\mu,\nu) := \inf_{\l\in J(\mu,\nu)}\l\{(\w,\eta) \in K^T\times K^T:\ \w_v \neq \eta_v\}\]
for any (arbitrary) choice of reference index $v \in T$, where
$J(\mu,\nu)$ denotes the collection of all \textbf{joinings} of
$\mu$ and $\nu$: $\G$-invariant probability measures on $K^T\times
K^T$ having first marginal $\mu$ and second marginal $\nu$.  If
$\rho(\mu,\nu)$ is small we shall write informally that $\mu$ and
$\nu$ have a \textbf{near-diagonal joining}.

The joining topology is clearly at least as strong as the vague
topology, and in general it is strictly stronger
(see~\cite{Ald82.2}, for example).  However, given a $\G$-invariant
closed subset $\O \subseteq K^T$, we can consider the subsimplex
$\Pr^\G\O \subseteq \Pr^\G K^T$ of exchangeable measures
concentrated on $\G$, and ask whether the two different
neighbourhood bases of this subsimplex defined by these two
topologies might coincide.  This question is motivated by the case
of hypergraph-exchangeability, for which it can be proved that these
bases \emph{do} always coincide; this follows, in particular, from
the rather more precise results for such closed subsets contained
in~\cite{AusTao--hereditarytest}. However, by making reference to
the Poulsen property, we can see that this is not always the case
for cube-exchangeability.

\begin{prop}\label{prop:Poulsenfromtestable}
If an exchangeability context $(T,\G)$ has the Poulsen property and
these two neighbourhood bases around $\Pr^\G\O$ are equivalent then
$\Pr^\G\O$ must also be the Poulsen simplex.
\end{prop}

\textbf{Proof}\quad In general, if $\mu_1$ is ergodic and is close
to $\mu_2$ in the vague topology, it need not follow that $\mu_1$ is
close to any of the ergodic components of $\mu_2$ in the vague
topology.  However, if in fact $\mu_1$ is joining-close to $\mu_2$
then it does follows that it is joining-close to many of the ergodic
components of $\mu_2$, by considering the ergodic decomposition of
the joining itself.

Let the situation be as described, and suppose that $\mu \in
\Pr^\G\O$; we must show that $\mu$ is vaguely approximable by
extreme points of $\Pr^\G\O$.  Since $\Pr^\G\O$ is just the subset
of those members of $\Pr^\G K^T$ that are concentrated on $\O$, its
extreme points are still just its ergodic members.

By the Poulsen property of $\Pr^\G K^T$, we know $\mu$ can be
vaguely approximated by ergodic measures in this larger simplex. On
the other hand, by the assumed equivalence of the two neighbourhood
bases, it follows that provided these approximating measures are
close enough to the subsimplex $\Pr^\G\O$ for the vague topology,
they actually have near-diagonal joinings with members of this
smaller simplex $\Pr^\G\O$.

However, if $\mu_1 \in \Pr^\G K^T$ is ergodic and $\l \in
\Pr^\G(K^T\times K^T)$ is a near-diagonal joining of $\mu_1$ to some
member of $\Pr^\G\O$, then the components of the ergodic
decomposition of $\l$ must (almost surely) be joinings of $\mu_1$ to
ergodic measures that are still members of $\Pr^\G\O$, and in order
that $\l$ be near-diagonal these ergodic components of $\l$ must
also be near-diagonal with high probability. It follows that $\mu_1$
must actually be joining-close, and hence vaguely close, to some
ergodic members of $\Pr^\G\O$; and since $\mu_1$ was itself vaguely
close to $\mu$, we deduce that $\mu$ must be vaguely approximable by
extreme points of $\Pr^\G\O$, as required. \qed

We suspect that the above implication cannot be reversed (in that
there are also $\O$ for which the neighbourhood bases do not
coincide, but for which $\Pr^\G\O$ is Poulsen anyway).

\begin{cor}
The subset $\O_2 \subseteq \bbF_2^{\bbF_2^{\oplus \bbN}}$ is such
that the joining neighbourhood basis of the simplex $\Pr^\G\O_2$ is
strictly stronger than the vague neighbourhood basis.
\end{cor}

\textbf{Proof}\quad By the previous proposition, it suffices to
argue that $\Pr^\G\O_2$ is not Poulsen; however, as discussed in the
previous subsection, the members of $\Pr^\G\O_2$ are precisely
Aldous' random walk examples in the case $U = \bbF_2$, and it is now
easy to check from this that the simplex in question has set of
extreme points precisely the measures $\mu$ constructed from $\nu_p$
for different $p > 0$ from the proof of
Lemma~\ref{lem:basicvagueapprox}, together with $\delta_\bf{0}$ and
$\delta_\bf{1}$, and that this set of extreme points has only the
one additional non-ergodic cluster point $\frac{1}{2}\delta_\bf{0} +
\frac{1}{2}\delta_\bf{1}$ (indeed, that lemma itself guarantees that
this must be cluster point; it is the argument of
Theorem~\ref{thm:Poulsen} that then necessarily takes us outside
$\Pr^\G\O_2$, and so does not apply to this sub-simplex). Thus,
$\Pr^\G\O_2$ cannot be Poulsen. \qed

In the setting of hypergraph exchangeability, it turns out that
there is a close relationship between properties of the sub-simplex
$\Pr^T\O$ and of the conditions on a point of $K^T$ needed to
guarantee membership of $\O$.  In addition, it turns out that this
latter membership condition can be identified simply with some
hereditary property of $K$-colourings of finite hypergraphs
(precisely, so that a point of $K^T$ lies in $\O$ if and only if
when regarded as a $K$-coloured hypergraph all of its finite induced
coloured sub-hypergraphs have that hereditary property).  From this
vantage point, a suitable analysis of this simplex can be converted
into a proof that all such properties are `efficiently testable'
(following essentially a translation of older, purely combinatorial
arguments to that effect; see, in particular, Alon and
Shapira~\cite{AloSha07} and R\"odl and Schacht~\cite{RodSch07}). We
shall not enter into these notions further here, but refer the
reader to the complete account in~\cite{AusTao--hereditarytest}.

It seems clear that a similar notion of efficient testability can be
formulated in the setting of discrete cubes and their isometries: in
general, we would write that a property $\P$ of all subsets of faces
of the finite discrete cubes $\bbF_2^N$ is \textbf{testable} if for
any $\eps > 0$ there are some $N(\eps) \geq J(\eps)\geq 1$ and
$\delta(\eps) > 0$ such that, if $N\geq N(\eps)$ and $E \subseteq
\bbF_2^N$, and if we know that a $J(\eps)$-face $F$ of $\bbF_2^N$
chosen uniformly at random has probability at least $1 -
\delta(\eps)$ of having $F \cap E \in \P$, then there is some $E'
\subseteq \bbF_2^N$ having $E' \in \P$ and $|E\Delta E'| < \eps
2^N$.

Although we are not aware of a rigorous relationship between the
question of Proposition~\ref{prop:Poulsenfromtestable} and
testability, by analogy with the results
of~\cite{AusTao--hereditarytest} we suspect from that Proposition
that the property $\O_2$ is \emph{not} testable; and in fact a
direct re-write of the particular infinitary proofs we have given in
finitary terms in a high-dimensional cube $\bbF_2^N$ shows that this
is so; we omit the details.

\subsection{Affine transformations of the infinite-dimensional discrete cube}

We have already discussed cube-exchangeability as a strengthening of
the condition of hypergraph-exchangeability treated by classical
exchangeability theory.  However, it may be worth recalling that an
even stronger exchangeability context on $T = \bbF_2^{\oplus\bbN}$
has also appeared implicitly in a number of recent works, with $\G$
the group of all affine transformations of $T$.

In particular, this setting closely relates to several questions of
current interest in arithmetic combinatorics concerning the counting
of affine copies of various patterns (such as finite-dimensional
cubes) in subsets of $\bbF_2^N$ for large $N$. These questions often
correspond naturally to descriptions of probability measures on
$\{0,1\}^{\bbF_2^{\oplus\bbN}}$ that are
$\rm{Aff}\,\bbF_2^{\oplus\bbN}$-invariant via a suitable
correspondence principle, analogous to the well-known Furstenberg
correspondence principle relating subsets of $\bbZ$ to
measure-preserving $\bbZ$-actions (see, for example, Furstenberg's
book~\cite{Fur77}). Closely-related to this line of research is the
investigation of the `Gowers-inverse conjecture' of Green and Tao in
the case of the vector spaces $\bbF_2^N$, which are phrased in terms
of correlations of individual $\bbC$-valued functions on $\bbF_2^N$
with functions of certain special forms. However, this conjecture
has recently been shown to fail in general in this setting in the
paper~\cite{GreTao07} of Green and Tao, and so some more complicated
kinds of ingredient seem to be required for such a structure
theorem.

In our more infinitary set-up, we suspect that in the presence of
this rather stronger symmetry a much more detailed analysis of the
structure of the exchangeable measures is possible, and that such an
analysis will probably rely on more ergodic-theoretic tools (such as
those developed for the proof or convergence and expression of the
limit of nonconventional ergodic averages in the case of
$\bbZ$-systems; see, in particular, the works of Host \&
Kra~\cite{HosKra05} and Ziegler~\cite{Zie04}); however, we have not
investigated this possibility further.  We also direct the reader to
Subsection 4.7 of~\cite{Aus--ERH} for a very informal discussion of
the different approaches to the extraction of structural information
for invariant measures in the study of exchangeability, on the one
hand, and ergodic theory on the other.

\subsection{The Poulsen property for other exchangeability contexts}

We suspect that the conclusion of Theorem~\ref{thm:Poulsen} holds
much more generally: that for an amenable group $\G$ it is only in
the presence of some very special exchangeability context (such as
those that are DMT) that the Poulsen property fails.

Is it possible to formulate a more general condition under which an
exchangeability context has the Poulsen property that will subsume
Theorem~\ref{thm:Poulsen}? On the other hand, is there some
condition related to that of being DMT that is actually equivalent
to the Bauer property (possibly only for amenable $\G$)? Can the
simplex $\Pr^\G K^T$ ever be neither Bauer nor Poulsen?

\subsection{Cube-exchangeability for finer-grained cubes}

We suspect that the results of this paper extend to the analogous
definition of exchangeability on the finer-grained cubes
$(\bbZ/m\bbZ)^{\oplus \bbN}$ for $m > 2$ (indeed, the situation
there is surely even more wild, if anything), but it is not clear
whether these exhibit any additional new phenomena.

\parskip 0pt

\bibliographystyle{abbrv}
\bibliography{Infinite_cube_exchangeability}

\vspace{10pt}

\textsc{Department of Mathematics\\ University of California at Los
Angeles\\ Los Angeles, CA 90095-1555, USA}

\vspace{7pt}

Email: \verb|timaustin@math.ucla.edu|

Web: \verb|http://www.math.ucla.edu/~timaustin|

\vspace{7pt}


\end{document}